  \documentclass[final,3p,times]{elsarticle}

\usepackage{amsfonts,amsmath,amsthm,amssymb}
\usepackage{epsfig}
\usepackage{verbatim}
\usepackage{color}

\def\NN{\hbox{I\kern-.2em\hbox{N}}}
\def\RR{{\mathop{{\rm I}\kern-.2em{\rm R}}\nolimits}}
\def\Q{{\bf Q}}
\def\C{{\bf C}}
\def\x{{\bf x}}
\def\y{{\bf y}}
\newcommand{\be}{\begin{equation}}
\newcommand{\ee}{\end{equation}}
\newcommand{\ba}{\begin{eqnarray}}
\newcommand{\ea}{\end{eqnarray}}

\begin{document}

\begin{frontmatter}



\title{Isogemetric Analysis and Symmetric Galerkin BEM:\\
a 2D numerical study}


\author[label1]{A.~Aimi\corref{cor1}}
\ead{alessandra.aimi@unipr.it}
\address[label1]{Department of Mathematics and Computer Science, University of Parma,\\
Parco Area delle Scienze, 53/A, Parma, Italy}
\cortext[cor1]{Corresponding Author}
\author[label1]{M.~Diligenti}
\ead{mauro.diligenti@unipr.it}
\author[label2]{M.~L.~Sampoli}
\ead{marialucia.sampoli@unisi.it}
\address[label2]{Department of Information Engineering and Mathematics, University of
Siena,\\
Via Roma 56, Siena, Italy }
\author[label3]{A.~Sestini}
\ead{alessandra.sestini@unifi.it}
\address[label3]{Department of Mathematics and Computer Science, University of
Florence,\\
Viale Morgagni 67, Firenze, Italy }

\begin{abstract}
\noindent Isogeometric approach applied to Boundary Element Methods
is an emerging research area (see e.g. \cite{Simpson}).
In this context, the aim of the present contribution is that of
investigating, from a numerical point of view, the Symmetric
Galerkin Boundary Element Method (SGBEM) devoted to the solution of 2D boundary value
problems for the Laplace equation, where the boundary and the
unknowns on it are both represented by B-splines (\cite{dB}). We
mainly compare this approach, which we call IGA-SGBEM, with a
curvilinear SGBEM (\cite{Aimi1}), which operates on any boundary
given by explicit parametric representation and where the approximate
solution is obtained using Lagrangian basis. Both techniques
are further compared with a standard (conventional) SGBEM approach
(\cite{Aimi}), where the boundary of the assigned problem is
approximated by linear elements and the
numerical solution is expressed in terms of Lagrangian basis.
Several examples will be presented and discussed, underlying
benefits and drawbacks of all the above-mentioned approaches.
\end{abstract}

\begin{keyword}
Isogeometric Analysis, B-splines, Symmetric Galerkin Boundary
Element Method



\end{keyword}

\end{frontmatter}


\section{Introduction}
\label{introduzione}
Boundary Element Methods (BEMs) have become
an important technique for solving linear elliptic partial
differential equations (PDEs) appearing in many relevant physical and
engineering applications (e.g. potential, acoustics, elastostatics,
 etc.; see \cite{Krishnasamy, Sirtori, Wendland}). By means of the fundamental solution of the
considered differential operator, a large class of both exterior
and interior elliptic Boundary Value Problems (BVPs) can be
formulated as a linear integral equation on the boundary of the
given domain. The numerical analysis of these methods for 2D and
3D problems is now well studied (\cite{Aimi, Aimi1, Aimi0, Carini,
Krishnasamy}). The BEM can offer substantial computational
advantages over other approximation techniques, such as finite
elements (FEM) or finite differences (FDM).
Moreover, in some applications, the physically relevant information is not the solution inside the domain but rather its trace or its normal derivative on the domain boundary: these latter can be obtained directly from the numerical solution of Boundary Integral Equations (BIEs), whereas boundary values recovered from FEM solutions are not so accurate.
However, in order to achieve an
efficient numerical implementation of general validity, a number
of issues have to be dealt with special attention. One of the most
important, for the practical application of the BEM analysis, is
the evaluation of singular integrals over boundary elements. It is
only in these last decades that engineers and applied
mathematicians have started employing finite part integrals to
formulate several 2D and 3D BVPs, particularly in the area of
applied mechanics, as singular and hypersingular
BIEs in the so-called symmetric
formulations (\cite{Krishnasamy, Sirtori}). Symmetric Galerkin
Boundary Element Method (SGBEM)
- see \cite{bonnet, gray} and references therein for rather complete surveys -
is nowadays recognized as a valid
alternative BEM technique for the solution of boundary value
problems, yielding final symmetric
discretization matrices which are suitable for the coupling with
FEM
(\cite{Wendland2, Costabel0}).
Since its origins (\cite{Sirtori_0}),
a large amount of literature results has been produced for what concerns stability and convergence properties of the method related to smooth or piecewise smooth boundaries and to piecewise polynomial systems for the approximation of the BIE solution (see e.g. \cite{Wendland, Wendland1, Rannacher, babuska, postell, suri});
further, great effort has been devoted to the development of efficient
schemes for the double integration of weakly singular, Cauchy
singular and hypersingular integrals over boundary elements
(\cite{Aimi1, bonnet, Kane, Holzer}).

On the other side, isogeometric analysis (IGA) is a new method for the numerical treatment
of problems governed by PDEs. In its
first formulation introduced in the literature (\cite{HCB}), the aim
was to overcome some difficulties arising in FEMs, proposing a viable alternative to standard,
polynomial-based, finite element analysis.
Actually, the key issue in IGA is to retain the
description of the domain where the PDE is defined as it is given
by a Computer Aided Design (CAD) system (i.e., in terms of B-splines or their rational
generalization, NURBS) instead of approximating it by a
triangular/polygonal mesh. Note that the most domains of interest in
engineering problems are exactly described in terms of B-splines
or NURBS. The term {\it isogeometric} is due to the fact that the
solution space for dependent variables is represented in terms of
the same functions which describe the geometry of the domain.
Thus, the isogeometric approach ensures an exact description of
the domain at any level, no matter how coarse is the
discretization of the problem. In addition, the mesh refinement is
highly simplified because it can be obtained by standard
\textit{knot-insertion} and/or \textit{degree-elevation} procedures (\cite{handbook}),
retaining the exact geometry of the original domain during the
process and eliminating the need to communicate with the CAD
system, once the initial mesh is constructed. In addition, the
easier manipulation of smooth elements provides an efficient tool
particularly well suited for high order equations. The above
mentioned facts motivate the wide interest received by this new
paradigm since the seminal paper \cite{HCB} (see for example
\cite{MPS'11, SMPS'12, VGJS'11} and references
quoted therein).

Very recently, literature on IGA has started dealing with
applications involving BEMs (see e.g. \cite{Bem_vienna'14, politis,
politis'14, Simpson}), even if the germinal idea was pushed forward
in \cite{Cabral_1, Cabral_2}.  This new approach has been mostly
compared with standard versions of the BEMs, based on a piecewise
polygonal approximation of the boundary of the problem domain,
obtaining, with no doubt, a remarkable superiority.

In this framework, we investigate, from a numerical point of view, the
so-called IGA-SGBEM, i.e. the SGBEM
devoted to the solution of 2D BVPs for the
Laplace equation (but the analysis could be extended considering
other operators), where B-splines are used to approximate the
boundary geometry as well as the unknown potential and flux fields
on it.

In order to combine IGA and SGBEM, we choose to work with
B--spline basis, since it is a fundamental tool for dealing with
polynomial spline spaces in the context of CAD and automatic
manufacturing, where spline functions expressed in the B--spline
basis (B--form)  are the standard for free--form design
(\cite{handbook}). Furthermore, B-splines are also profitably used
in several other fields,  for example in multiresolution analysis
or in collocation methods (see e.g. \cite{LMQ01, MST06}) and, as
above mentioned, they are also the standard basis adopted in the
recent context of IGA-FEM (\cite{IgAbook}). The reason of this wide
success depends on several aspects, surely because splines in
B--form can be easily stored, evaluated and algebraically
manipulated.
Another important reason is that they
have general features, very attractive for applications, such as
nonnegativity, partition of unity, compact support, and
total positivity (\cite{dB, CP94}).

In this paper, IGA-SGBEM approach will be compared with
curvilinear SGBEM (\cite{Aimi1}) - an improved version with respect
to the existing conventional  boundary element techniques - which
operates on any boundary given by explicit parametric representation
(hence, in particular, given by
B-splines representation), and where the approximate
solution is obtained using Lagrangian basis. Both the above
mentioned methods will be further compared with a standard
(conventional) SGBEM approach (\cite{Aimi}), where the boundary of
the assigned problem is approximated by linear elements and the
numerical solution is expressed in terms of Lagrangian basis.
Singular integrals required by SGBEMs are efficiently evaluated by
suitable quadrature formulas with a very low number of quadrature
nodes related to user assigned accuracy (\cite{Aimi, Aimi1}).
Several examples will be presented and discussed, underlying
benefits and drawbacks of all the above mentioned techniques.

\section{Model problem and its boundary integral formulation}
\label{} Let $\Omega \subset \mathbb{R}^2$ be a bounded, simply
connected, open domain with a (piecewise) smooth boundary $\Gamma:=\partial \Omega=\{\x=(x_1,x_2) \in
 \mathbb{R}^2 | \,\x= \C(t),\, t\in [a,b]\}$, given by parametric representation on the interval $[a,b]$.
 Let us further suppose that
$\Gamma={\bar \Gamma}_1\cup{\bar \Gamma}_2$, where $\Gamma_1$ and $\Gamma_2$ are
open disjoint subset of $\Gamma$
and $\textit{meas}\,(\Gamma_1)>0$.
As model problem, we consider a
mixed BVP for the Laplace equation:

\indent \textit{given $u^* \in H^{1/2}(\Gamma_1)$ and $q^*\in H^{-1/2}(\Gamma_2)$,
find $u\in H^1(\Omega)$ such that}
\be \label{uno}
\left\{
\begin{array}{ll}
\Delta u=0 & \textrm{in} \,\,\Omega\,,\\
u=u^* & \textrm{on} \,\,\,\Gamma_1\,,\\
\frac{\partial u}{\partial \mathbf{n}}=q^* & \textrm{on}
\,\,\,\Gamma_2\,,
\end{array}
\right. \ee
where $\frac{\partial}{\partial \mathbf{n}}$ denotes the
derivative with respect to the outer normal $\mathbf{n}$ to
$\Gamma$. The definition of Sobolev spaces is as usual (see
\cite{Lions}). \\
As it is well known (\cite{Aimi0, Costabel, Wendland}), from problem $(\ref{uno})$ the following identities for
$u$ and $q$ on $\Gamma$ can be derived:
\begin{equation}\label{quattro}
\frac{1}{2}\left[\begin{array}{c}u \\ q \end{array}\right]=\left[\begin{array}{cc}-K & V \\
-D & K'\end{array}\right]\left[\begin{array}{c}u \\ q
\end{array}\right]\,, \qquad \x\in\Gamma\,,
\end{equation}
where
$$
\begin{array}{ll}
\displaystyle
Vq(\x):=\int_{\Gamma}U(\x,\y)\,q(\y)\textrm{d}\gamma_\y\,,&
Ku(\x):=\displaystyle\int_{\Gamma}
\frac{\partial U}{\partial \mathbf{n}_\y}\,(\x,\y)\,u(\y)\textrm{d}\gamma_\y\\
\\
K'q(\x):=\displaystyle\int_{\Gamma}
\,\frac{\partial U}{\partial
\mathbf{n}_\x}\,(\x,\y)\,q(\y)\textrm{d}\gamma_\y\,,&
Du(\x):=\displaystyle\int_{\Gamma}
\,\frac{\partial^2 U}{\partial \mathbf{n}_\x\partial
\mathbf{n}_\y}\,(\x,\y)\,u(\y)\textrm{d}\gamma_\y\,,
\end{array}
$$
and
$$
U(\x,\y):=-\frac{1}{2\pi}\ln \|\y-\x\|_2\,,
$$
is the fundamental solution of the 2D Laplace operator. Note that
$K$ and $K'$ are defined by Cauchy singular integrals when $\Gamma$ is a piecewise smooth boundary \footnote{In the case of a smooth boundary, the operators $K$ and $K'$ are only weakly singular (\cite{Attkinson}, Section 7).}, while $D$ is
defined by a hypersingular finite part integral in the sense of
Hadamard, i.e. it is understood to be the finite part of an asymptotic expansion (\cite{Schwab,
 Wendland}).\\ Under the above assumptions, the
following properties are well known (\cite{Costabel,
 Schwab, Wendland}): the operators
\begin{equation}
\begin{array}{ll}
V:H^{-1/2+\sigma}(\Gamma)\,\rightarrow\,H^{1/2+\sigma}(\Gamma)\,,
\quad & K:H^{1/2+\sigma}(\Gamma)\,\rightarrow\,H^{1/2+\sigma}(\Gamma)\,,\\
K':H^{-1/2+\sigma}(\Gamma)\,\rightarrow\,H^{-1/2+\sigma}(\Gamma)\,,\quad&
D:H^{1/2+\sigma}(\Gamma)\,\rightarrow\,H^{-1/2+\sigma}(\Gamma)\,,
\end{array}
\label{opint}
\end{equation}
are continuous for $\sigma \in (-\frac{1}{2},\frac{1}{2})$. For
$\sigma=0$ the operator $K'$ is the adjoint of $K$ with respect to
the natural duality $<\cdot,\cdot>$ between $H^{1/2}(\Gamma)$ and
its dual $H^{-1/2}(\Gamma)$,
which for sufficiently smooth functions coincides with the usual scalar product in $L^2(\Gamma)$.

The strong system of two BIEs $(\ref{quattro})$  is, of course, overdetermined (see e.g. \cite{bonnet, Wendland}), hence it can be
reformulated without redundancy, following an approach similar to the one rigorously developed in
\cite{Andra}. Therefore, by imposing the first equation
only on $\Gamma_1$ and the second one only on $\Gamma_2$, and inserting the boundary data given in
$(\ref{uno}),$ we obtain a system of two BIEs of the first kind
for the unknowns $q$ on $\Gamma_1$ and $u$ on
$\Gamma_2$, of the form
\begin{equation}\label{cinque}
\left[\begin{array}{cc} V_{11} & -K_{12} \\ -K'_{21} &
D_{22}\end{array}\right] \left[\begin{array}{c}q \\
u\end{array}\right]=\left[\begin{array}{cc} -V_{12} &
\frac{1}{2}I+K_{11} \\ -\frac{1}{2}I+K'_{22} &
-D_{21}\end{array}\right] \left[\begin{array}{c}q^* \\
u^*\end{array}\right]\,,
\end{equation}
where the boundary integral operators subscripts \textit{j\,k}
mean evaluation over $\Gamma_j$ and integration over $\Gamma_k$. Note that an alternative non redundant system
of two BIEs of the second kind could be obtained considering the first equation only on
$\Gamma_2$ and the second one only on $\Gamma_1$. However this latter
approach does not lead to symmetric final discretization matrices.\\
System $(\ref{cinque})$ will be solved in a weak sense (see
\cite{Aimi0, Wendland}). The weak formulation starts from
identity $(\ref{quattro})$:  finding the weak solution $u\in
H^1(\Omega)$ of BVP $(\ref{uno})$ is indeed equivalent to find the
weak solution $[u,\,q] \in H^{1/2}(\Gamma)\times H^{-1/2}(\Gamma)$
of system $(\ref{quattro})$ such that $u|_{\Gamma_1}=u^*$ and
$q|_{\Gamma_2}=q^*$. \\ After having recovered the
missing Cauchy data by solving, with obvious meaning of notation,
the weak symmetric (\cite{bonnet}) problem:
\begin{equation}\label{cinquebis}
<\left[\begin{array}{cc} V_{11} & -K_{12} \\ -K'_{21} &
D_{22}\end{array}\right] \left[\begin{array}{c}q \\
u\end{array}\right],\left[\begin{array}{c}p \\
v\end{array}\right]>=< \left[\begin{array}{c}f_1 \\
f_2\end{array}\right],\left[\begin{array}{c}p \\
v\end{array}\right]>\,,\quad \quad \forall\, [p,v]\in
H^{-1/2}(\Gamma_1)\times H^{1/2}_0(\Gamma_2)\,,
\end{equation}
one can use the representation formula
$$
u(\x)=\int_\Gamma U(\x,\y)\,q(\y) \textrm{d}\gamma_\y -\int_\Gamma
\frac{\partial U}{\partial \mathbf{n}_\y}\,(\x,\y) \,u(\y)
\textrm{d}\gamma_\y\,,\qquad \x\in \Omega\,, $$
to obtain the solution at any point of the domain.\\

\noindent Let us remember the advantages of using the symmetric
boundary integral problem \eqref{cinque} for a mixed BVP: the
unknowns on the boundary are directly those of the differential
problem instead of density functions typical of indirect
formulations (\cite{chen}); the
linear system coming from the discretization of \eqref{cinquebis},
due to the involved integral operators properties, presents a
symmetric matrix and this, of course, is extremely important for
what concerns saving computational time in matrix
generation and memory in matrix storage.\\

\noindent \textit{Remark.} If we have to deal with a Dirichlet
BVP, i.e. $\Gamma \equiv \Gamma_1$, the systems \eqref{cinque}
obviously reduces to the first equation alone, where the only
unknown is $q(\x)$. A similar boundary integral equation can be
written for a Dirichlet problem exterior to an open arc in the
plane (see e.g. \cite{chen}): in this case, the unknown is the
jump of $q(\x)$ across the arc $\Gamma$, i.e. $[q({\bf x})]_\Gamma$.

\section{Symmetric Galerkin Boundary Element Method}
For the discretization phase, we consider a uniform partition of
the parametrization interval $[a,b]=\bigcup_{\ell=1}^n I_\ell$,
made up by $n$ subintervals $I_\ell$ and governed by the
decomposition parameter $h=length(I_\ell)$. This induces over
$\Gamma$, using the parametric representation of the boundary, a
mesh $\Gamma_h=\bigcup_{\ell=1}^n e_\ell$, constituted by
curvilinear elements $e_\ell=\C(I_\ell)$.\\  In a similar way, a
finite dimensional subspace of piecewise polynomial functions can
then be lifted on the boundary, starting from the introduced
partition of $[a,b]$.\\ In the IGA-SGBEM the very same B-spline
basis used to represent the boundary $\Gamma$ is used also as a
basis for the functional approximation space; in the curvilinear
SGBEM, the boundary can be given by any explicitly defined
parametric representation (and therefore also by B-spline
representation), but the approximation space is spanned by a
Lagrangian basis defined over the decomposition of $[a,b]$.\\ At
last, in the standard SGBEM, $\Gamma$ is approximated by a
polygonal boundary ${\tilde \Gamma}_h$, constituted by linear
elements, each interpolating the endpoints of $e_\ell,\, \ell=1,
\cdots, n$, and a local Lagrangian basis is lifted onto each
straight element of ${\tilde \Gamma}_h$ from the reference element
$[0,1]$.

In any case, denoting with $\{\phi_i\}$ the basis of the
functional approximation space where we will search the unknowns
by means of Galerkin criteria, the elements of the final
discretization linear system matrix will be double integrals of
the form
\begin{equation}\label{nucleo}
\int_{{\Gamma}_h} \phi_j(\x)\int_{{\Gamma}_h}
{\cal K}(\x,\y)\phi_i(\y) \, d\gamma_\y d\gamma_\x\,,
\end{equation}
(substituting ${\Gamma}_h$ with ${\tilde \Gamma}_h$ in the case of
standard SGBEM), where ${\cal K}$ denotes one of the kernels of
the integral operators \eqref{opint} and therefore it can be
weakly singular, singular or hypersingular. Consequently, the
inner integral in \eqref{nucleo} has to be defined as generalized,
Cauchy principal value or Hadamard finite part, respectively.\\
Then, using suitable numerical integration schemes to face all
these types of singularities (see \cite{Aimi, Aimi1}), one can
write down the linear, symmetric, non singular system of equations
\begin{equation}\label{cinqueter}
\left[\begin{array}{cc} R_{11} & R_{12} \\ R_{21} &
R_{22}\end{array}\right] \left[\begin{array}{c}q_h \\
u_h\end{array}\right]= \left[\begin{array}{c}b_1 \\
b_2\end{array}\right]\,,
\end{equation}
where the vector unknowns $q_h, u_h$ collect the coefficients with
respect to the selected basis, which allow to finally obtain an
approximate solution of the integral problem.\\

\vspace{-0.1in}
In the remaining part of the Section, we recall the definition of
B-spline basis used in IGA-SGBEM approach, because it could be not
necessarily known to people in the numerical simulation community.
On the opposite, we do not  recall the definition of the
Lagrangian basis used below to compare IGA-SGBEM with curvilinear
and standard SGBEMs, because it's more basic (anyway the
interested reader can refer to \cite{brebbia}).

Given a partition $\Delta:=\{a=t_0<\cdots<t_n=b\}$ of an interval $[a,b]$,
a general polynomial spline space
$S$ of order $k$ on such partition is composed by piecewise
polynomial functions of degree $k-1$ which are required to have an
assigned regularity $C^{k-1-m_i}$ at the breakpoints $t_i,
i=1,\ldots,n-1,$ with $m_i$ denoting an integer between $1$ and
$k$ \footnote{When $m_i=k$ this means that the function has a
finite jump at $t = t_i.$ }.
For example when all the $m_i$ are fixed equal to $1$ or to $k-1$
or  to $k,$ respectively    $S$ is a subset of $C^{k-2}[a,b]$, it
is included in $C[a,b]$ or it is just a subset of $L^2(a,b).$ It
is quite easy to verify that the dimension of such space is
$
dim(S) \,=\,  k \,+\,{\sum_{i=1}^{n-1}} m_i\,.
$
\noindent The easiest way to define in $S$ a B-spline basis
$B_{i,k}(t), i=0,\ldots,N,$    with $N+1 = dim(S),$ is based on
the usage of a recursion formula and can be described through two
easy steps (\cite{dB}). The first step consists in associating to
$S$ an \textit{extended knot vector} $T=\{\tau_0,\cdots,\tau_{N+k}
\}$ whose elements constitute a non decreasing sequence of
abscissas, where $\{\tau_{k-1} ,\cdots,\tau_{N+1} \}$ are the {\it
internal knots} with $\tau_{k-1}= t_0,\, \tau_{N+1}= t_n$ and
$\{\tau_k ,\cdots,\tau_N \} = \{
t_1,\cdots,t_1,\cdots,t_{n-1},\cdots,t_{n-1} \},$ where each $t_i$
has $m_i$ occurrences and it is said {\it multiple} if $m_i>1.$
The remaining knots in $T,$  $\{\tau_0, \cdots,\tau_{k-2} \}$ and
$\{\tau_{N+2}, \cdots,\tau_{N+k} \}$ form two sets of $(k-1)$
knots called {\it auxiliary left} and {\it right knots}  which are only
required to verify the following inequalities, $\tau_0 \le \cdots
\le \tau_{k-2} \le \tau_{k-1} =a$ and $b=\tau_{N+1} \le \tau_{N+2}
\le \cdots \le \tau_{N+k}\,.$ Note that, in the numerical
simulations, we will always use the standard assumption of
selecting an {\it open} extended knot vector, that is
$\tau_0=\cdots=\tau_{k-2} = \tau_{k-1} = a$ and  $b= \tau_{N+1} =
\tau_{N+2} = \cdots = \tau_{N+k}\,.$

\noindent In the second step, the basis
is defined by using the following recursion (\cite{dB}):
$$
 \begin{array}{ll}
B_{i,1}(t) &:=\, \left\{\begin{array}{ll}
1, & \mbox{if }  \tau_i \leq t< \tau_{i+1}\,,\\
0, & \mbox{otherwise.}
\end{array}\right.
\cr & \cr B_{i,j}(t) &:=\, \omega_{i,j}(t)\,B_{i,j-1}(t) + (\,1-
\omega_{i+1,j}(t)\, )\, B_{i+1,j-1}(t)\,, \quad 1 < j \le k\,,
 \end{array} $$
with $ \omega_{i,j}(t) \,:=\,  \left\{\begin{array}{ll}
\frac{t-\tau_i}{\tau_{i+j-1} - \tau_i}
& \mbox{ if }  \tau_i < \tau_{i+j-1} \,, \cr 0 & \mbox{ otherwise.} \end{array} \right.$\\
Note that from the above recursive definition it is easy to verify
the nonnegativity of B-splines and that the support of $B_{i,k}$
is the subinterval $[\tau_i\,,\,\tau_{i+k}].$ The partition of
unity property can also be easily proved by induction on the
order.\\

As an example, in Figure \ref{basi}, the plots of all B-splines
spanning  two different  quadratic spline spaces $S_1$ and $S_2$,
respectively of dimension  $N+1 = 13$  and $N+1 = 15$, can be seen.
Such spaces share the same partition  $\{t_i = i, i=0,\ldots,9 \}$
of the interval $[0\,,\,9]$ but their extended knot vectors   are $T_1$ on
the left and $T_2$ on the right, with
\begin{equation} \label{T1T2}
\begin{array}{ll}
T_1 &=\, [0,\, 0,\, 0,\, 1,\, 1,\, 2,\, 3,\, 4,\, 5,\, 6,\, 7,\, 8,\,
8,\, 9,\,9,\, 9] \,, \cr
T_2 &=\, [0,\, 0,\, 0,\, 1,\, 1,\, 1,\, 2,\, 3,\, 4,\,5,\,6,
7,\,8,\,8,\,8,\,9,\,9,\,9]\,. \cr
\end{array}
\end{equation}


\section{Numerical results}

From now on, we will indicate curvilinear and standard SGBEMs with C-SGBEM and
S-SGBEM, respectively.\\

\noindent {\bf Example 1.}

\begin{figure}[bht]
\centerline{\includegraphics[width=2.5in]{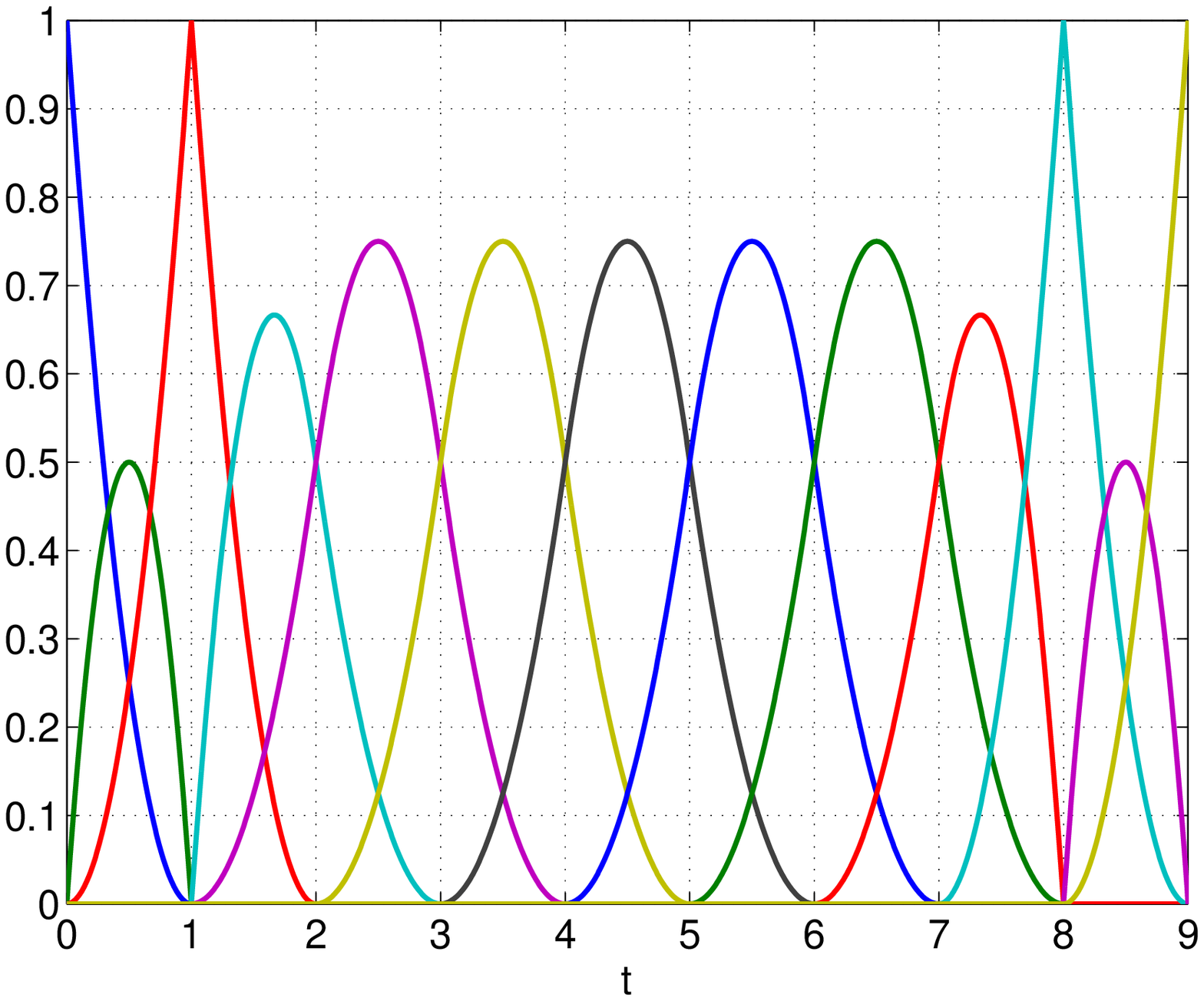}
\hskip.5cm\includegraphics[width=2.5in]{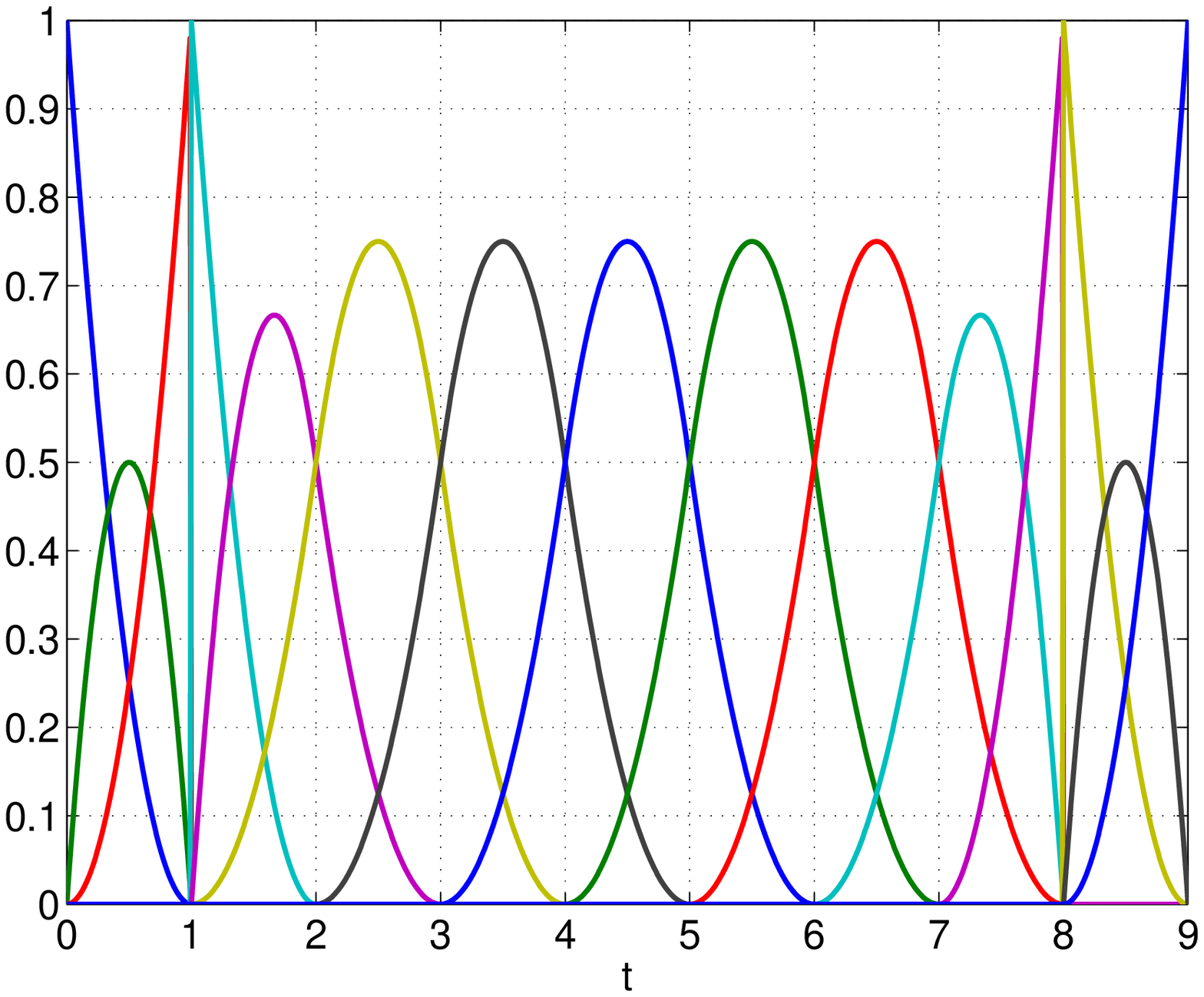}}
\caption{Quadratic B-splines basis with knot vector $T_1$ (left),
and knot vector $T_2$ (right), related to the same partition of the interval $[0,9]$.} \label{basi}
\end{figure}

\begin{figure}
\centerline{\includegraphics[width=3.in]{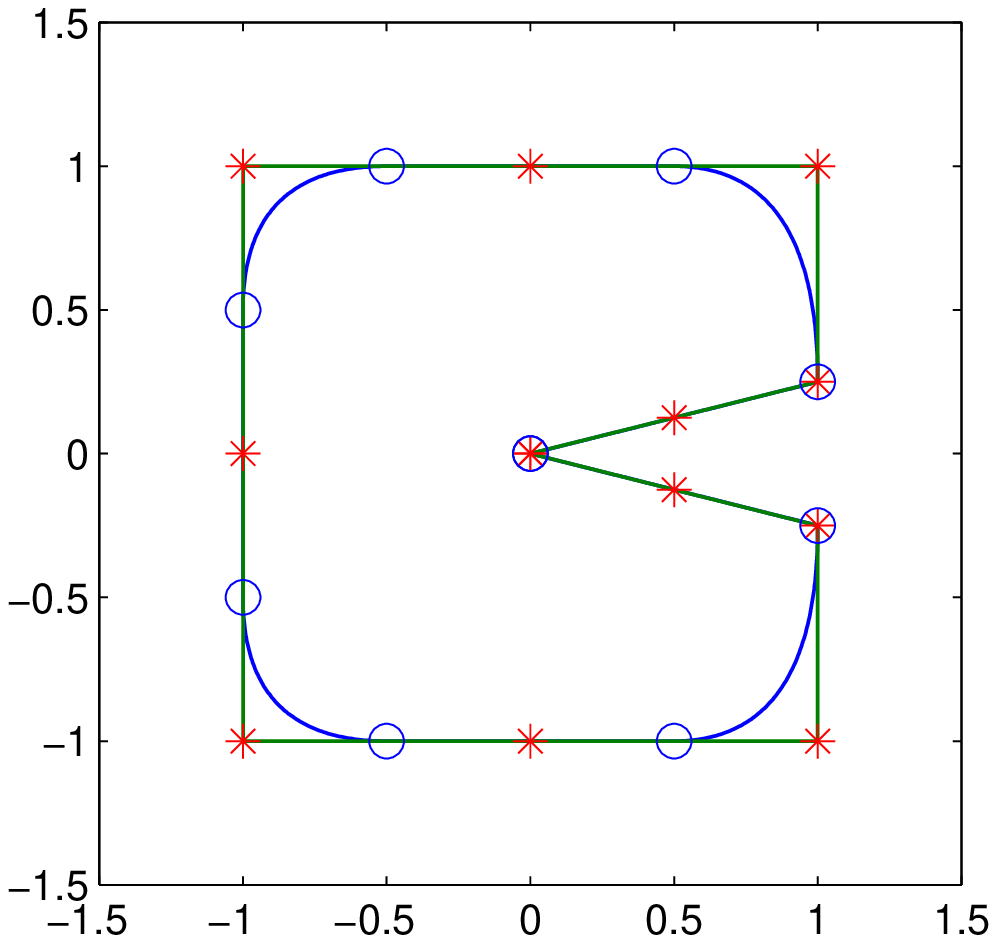}}
\caption{Example 1: the quadratic closed spline curve defining the
boundary of the considered interior domain along with the related
B-spline control polygon. The control points and the nodal (mesh) points
of the spline curve are respectively marked with the symbol
'$\ast$' and '$\circ$'. } \label{dominio}
\end{figure}

In this example we consider a potential problem interior to the
domain shown in Figure \ref{dominio}, which has three sharp
corners and is similar to that one constructed by means of NURBS
in \cite{Simpson}. The boundary of the domain is in our case
described by a closed parametric piecewise quadratic curve,
$\C(t),\, t \in [0,9]$, with integer uniform breakpoints $t_i =
i,\, i=0,\ldots,9.$ Such curve can be represented in B-form, that
is as a linear combination of the quadratic B-splines (see Figure
\ref{basi}, left) associated to the extended knot vector $T_1$
given in (\ref{T1T2}) and control points $\Q_i, \, i=0\ldots, 12$,
whose coordinates are given in the following matrix:
$$
Q=\Big[\begin{array}{c c c c c r r r r r c c c}
0  & 0.5&     1&   1&   0&  -1&  -1&  -1&  0&  1&   1 &  0.5&  0\\
0 &  0.125&     0.25&   1 &  1 &  1 &  0 & -1 & -1 &  -1&  -0.25&
-0.125 & 0
\end{array}\Big]\,.
$$
Note that the curve is closed because $\Q_0=\Q_{12}$ but its regularity at the initial/final
joint point is only $C^0$ because $T_1$ is an open extended knot vector.
Moreover, considering the double multiplicity of the breakpoints $t_1$ and $t_8$ specified
in $T_1,$ it turns out that  $\C(t) \in C^0[0,9]\cap C^1[1,8]$ .
In this way, the geometry of the domain boundary can be exactly described.

The differential problem is equipped with Dirichlet boundary
condition $u^*(\x)=-(x_1+x_2)$; the solution of the related
boundary integral equation is explicitly known, it reads
$q(\x)=q(\C(t))=(C'_1(t)-C'_2(t))/\|\C'(t)\|_2$ and has $L^2(0,9)$
regularity. In particular it presents a jump discontinuity at $t=
t_1,\, t=  t_8$ and at $t=t_0$ ($t=t_9$), while it is only $C^0$
in the remaining breakpoints.

As a first choice, we do not care about the low regularity of the
solution and we work in the space used to describe the boundary
which is spanned by the quadratic B-splines associated to $T_1$
and is a subset of $C^0[0,9]\cap C^1[1,8]. $  Then we successively
extend the space by inserting a new simple knot at the midpoint
between any two successive breakpoints (this corresponds to
halving the mesh step $h$, since uniform distributions of the
breakpoints are always assumed).

\noindent In Table \ref{ex1_tab1}, the obtained results are shown:
for each considered $h,$  the corresponding total number of
degrees of freedom (DoF), the spectral condition number of the
matrix in (\ref{cinqueter}), and the relative error \be
E=:\|q-q_h\|_{L^2}/\|q\|_{L^2}\,,\label{E} \ee are given.

\begin{figure}[hbt]
\centerline{\includegraphics[width=3.in]{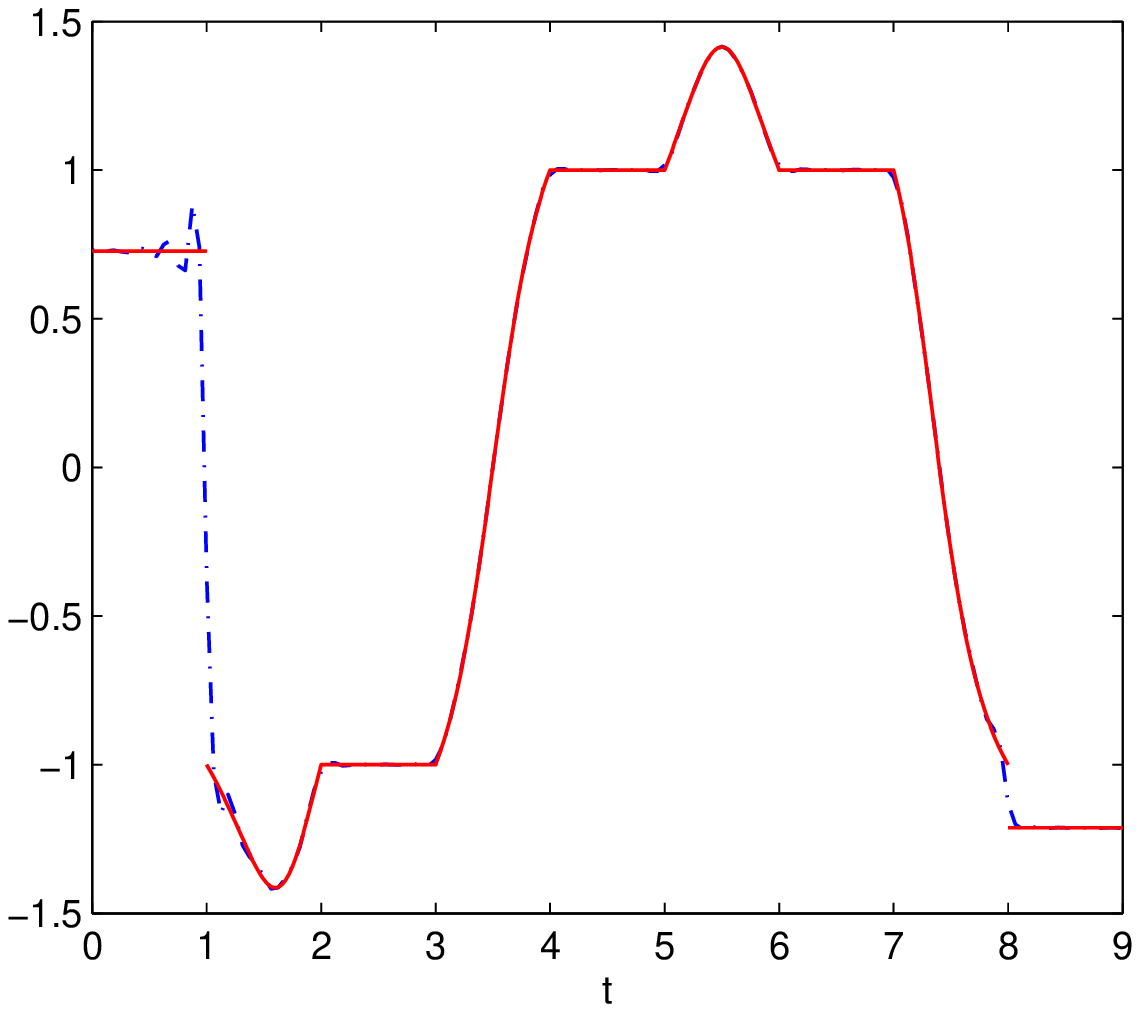}}
\caption{Example 1: the analytical solution (solid) and the
numerical solution (dash-dotted), obtained after three refinements of
$T_1$ ($h=1/8$).} \label{flusso_packman_quadratiche_ref3_C0_C1}
\end{figure}

\noindent Figure \ref{flusso_packman_quadratiche_ref3_C0_C1} confirms that
the numerical solution obtained with $h=1/8$ mainly agrees with
the analytical solution. As expected, the jumps are smoothly
approximated; small oscillations occur in the
neighborhood of the jumps, especially around $t=1$.

In order to adequate the quadratic spline space to the regularity
of the analytical solution,  we have then performed a similar set
of experiments starting now from the extended knot vector $T_2$,
given in (\ref{T1T2}).
The associated $L^2(0,9)\cap C^1[1,8]$ quadratic B-spline basis is
shown in Figure~\ref{basi}, right.

\noindent In Table \ref{ex1_tab2}, we show the comparison between the
results obtained successively refining the parameter $h$ for the $L^2(0,9)\cap
C^1[1,8]$ quadratic B-splines in IGA-SGBEM and the $L^2(0,9)\cap C^0[1,8]$ quadratic
Lagrangian basis in C-SGBEM. In
particular, for both approaches, we present the total number of
degrees of freedom (DoF), the spectral condition number of the associated
linear system matrix and the relative error \eqref{E}.

\begin{figure}[htb]
\centerline{\includegraphics[width=3.in]{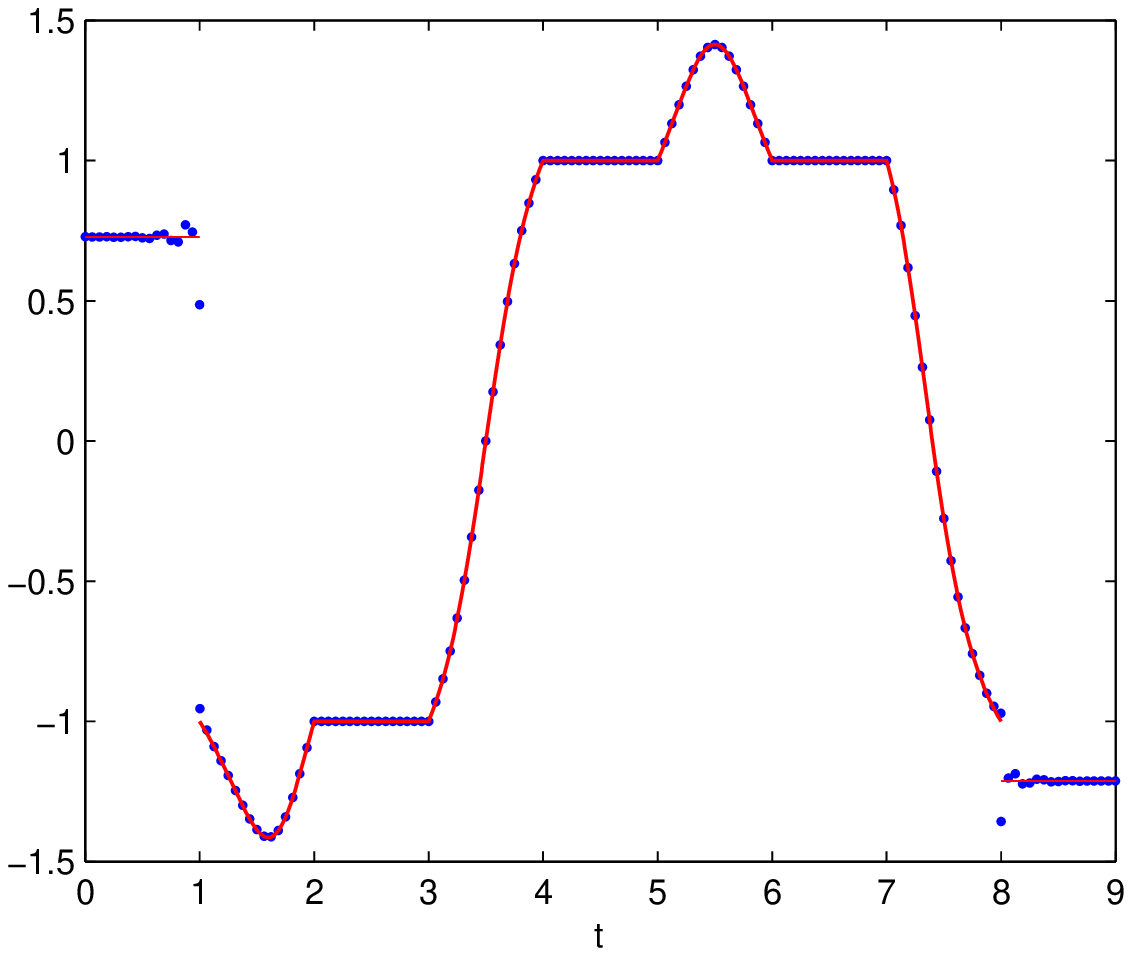}}
\caption{Example 1: the analytical solution (solid) and the
numerical solution (dotted), obtained after three refinements of $T_2$ ($h=1/8$).}
\label{flusso_packman_quadratiche_basic_ref3}
\end{figure}

\noindent Final errors are better than the corresponding ones shown in
Table \ref{ex1_tab1}, even if we note an error stagnation, for either IGA-SGBEM and C-SGBEM. This is due to
the difficulty of recovering the analytical solution
near the jumps, in particular on the side of the jumps where the exact solution is constant, even if, anyway, the approximation
sensibly improves elsewhere, as it is shown in Fig.~
\ref{flusso_packman_quadratiche_basic_ref3}. Actually, for this kind of solutions, as
well as for solutions exhibiting sharp layers,
the use of generalized exponential spline
spaces could be more suitable (\cite{MPS'11}) and it is planned as future work.

\medskip

\noindent {\bf Example 2.}

\begin{figure}[hbt]
\centerline{\includegraphics[width=3.in]{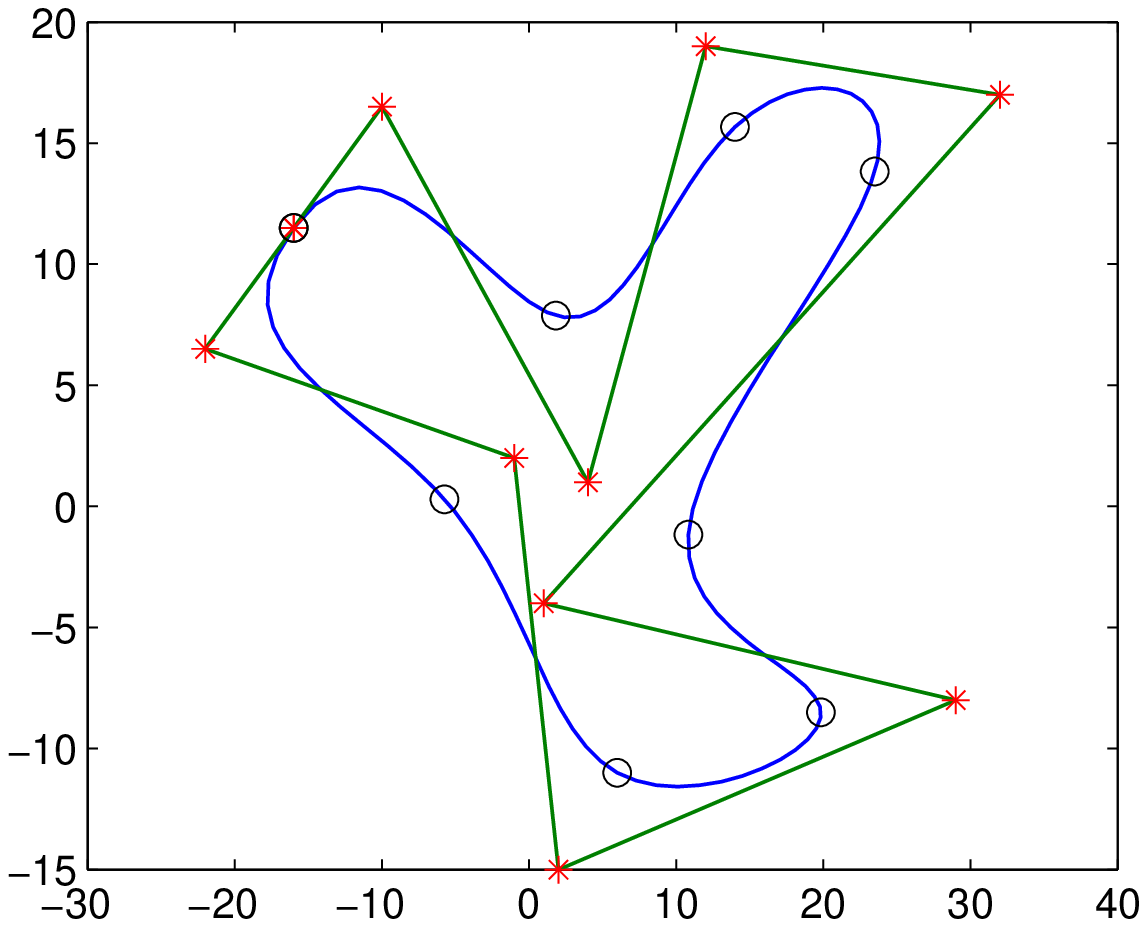}}
\caption{Example 2: the cubic closed spline curve defining the
boundary of the considered interior domain along with the related
B-spline control polygon. The control points and the nodal (mesh) points
of the spline curve are respectively marked with the symbol
'$\ast$' and '$\circ$'.} \label{free_form}
\end{figure}

In the second example we consider a potential problem
interior to the domain shown in Figure \ref{free_form}, see
\cite{politis}. Such domain has a smooth boundary that can be
described by a closed parametric piecewise cubic curve, $\C(t),\,
t \in [0,1]$, with uniform breakpoints and mesh step $h=1/8$.
This curve can be represented in B-form with extended knot
vector
$$
T_3= [0,\,0,\, 0,\, 0, \, 0.125,\, 0.25,\, 0.375,\, 0.5,
0.625,\,0.75,\,0.875,\,1,\,1,\,1,\,1 ]$$ and control points $\Q_i,
\, i=0\ldots, 10$, whose coordinates are collected in the following matrix:
$$
Q=\Big[\begin{array}{c r r r r r r c c c c c c c c}
-16  & -22&   -1&   2&   29&   1&  32&  12&    4&  -10&  -16  \\
11.5 &  6.5&   2&   -15&   -8 &  -4 &  17 &  19 & 1  &  16.5&  11.5
\end{array}\Big].
$$

The differential problem is equipped with Dirichlet boundary
condition $u^*(\x)=-(x_1+x_2)$; as in the previous example, the
solution $q(\x),\,\x=\C(t)$, of the related boundary integral
equation is explicitly known, but now, as function of t, it is
$C^1$ regular on $[0,1]$.

For this example, we first present a comparison between the
results obtained working in nested $C^2$ spline spaces of
increasing degree $\ge 3$, spanned by the B-spline basis
(IGA-SGBEM), and working with larger $C^0$ spline spaces of
corresponding degree spanned by the Lagrangian basis (C-SGBEM).
Note that the boundary curve can be exactly expressed in all the
considered spaces and its representation can be obtained combining
a \textit{degree elevation} with a \textit{knot
insertion} procedures (see e.g. \cite{handbook}).

\begin{figure}
\centerline{\includegraphics[width=2.8in]{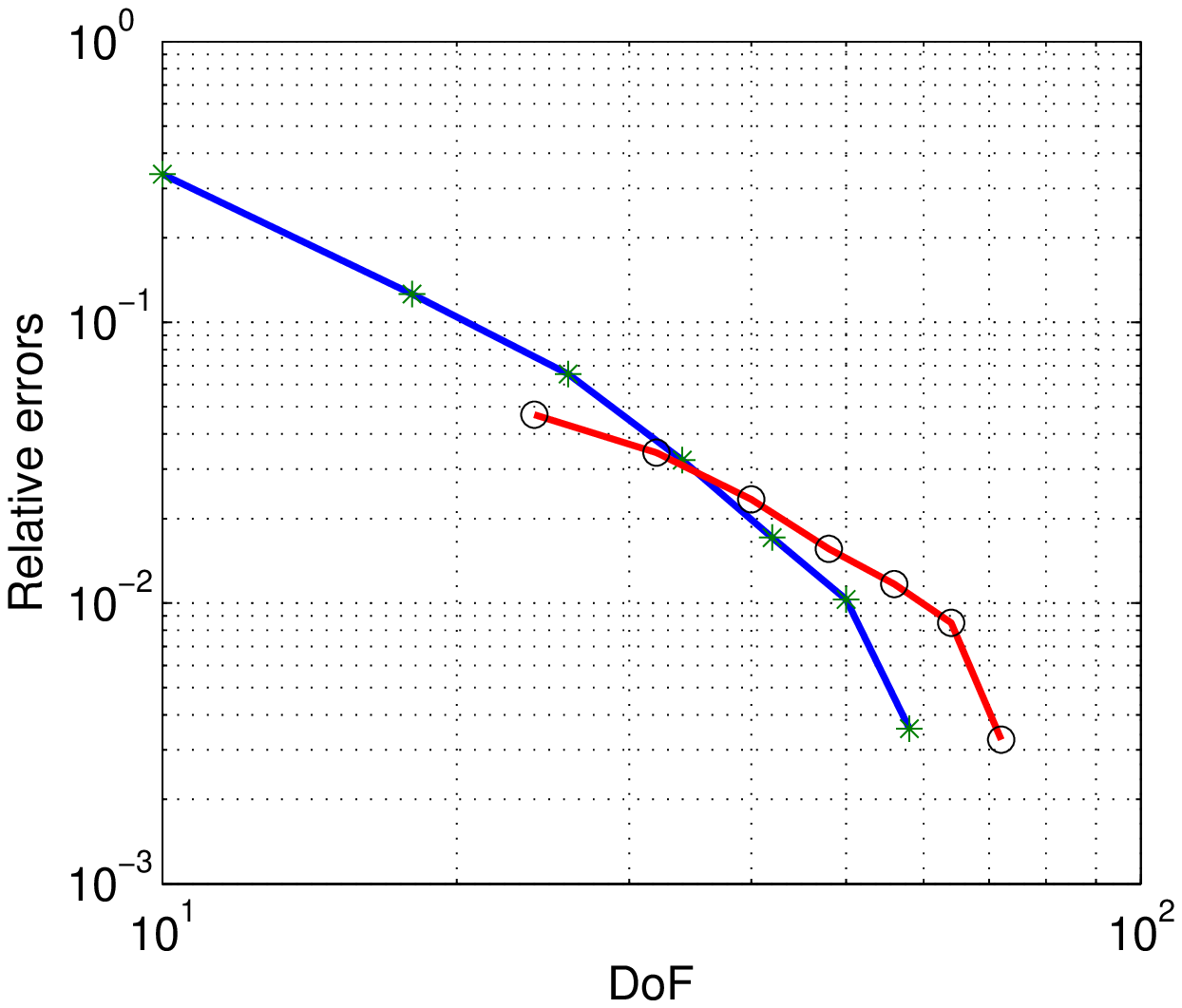}}
\caption{Example 2: relative errors of Table \ref{ex2_tab1}  vs
DoF with $h=1/8$ ('$*$' B-splines, '$\circ$' Lagrangian basis). }
\label{errorivsdof_basic}
\end{figure}

The recalled basis functions regularity, which
defines the type of norms that can be used to estimate the
approximation error (\cite{Rannacher}), allows us to compare the
different approaches considering the relative error $E$ defined
in \eqref{E}, as done in the previous example.
Results are presented in Table \ref{ex2_tab1} and the
errors are then plotted in Figure \ref{errorivsdof_basic} with
respect to DoF.

\begin{figure}
\centerline{\includegraphics[width=2.8in]{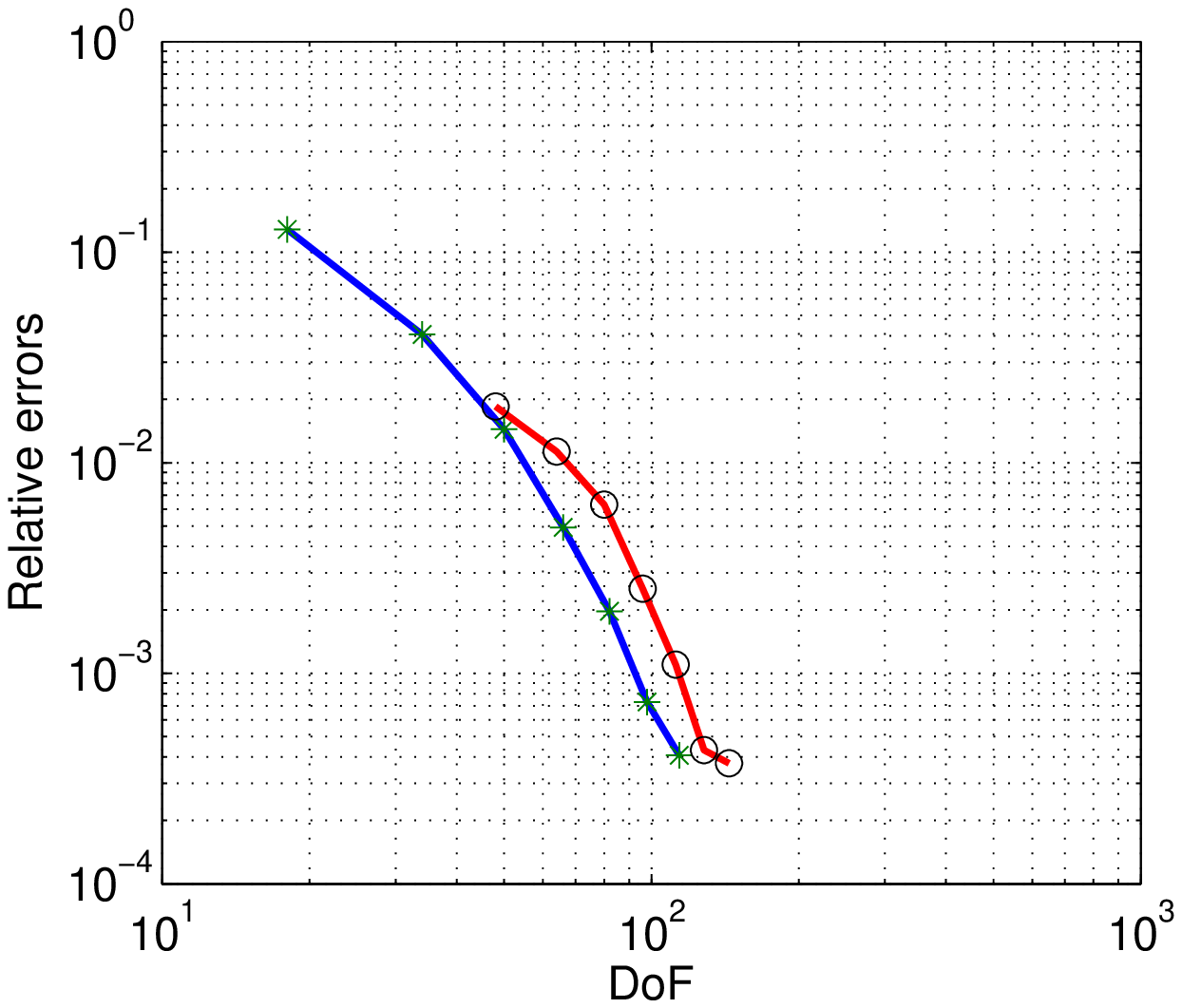}}
\caption{Example2: relative errors of Table \ref{ex2_tab2}  vs DoF
with $h=1/16$ ('$*$' B-splines, '$\circ$' Lagrangian basis).}
\label{errorivsdof_Iraff}
\end{figure}

\begin{figure}
\centerline{\includegraphics[width=2.8in]{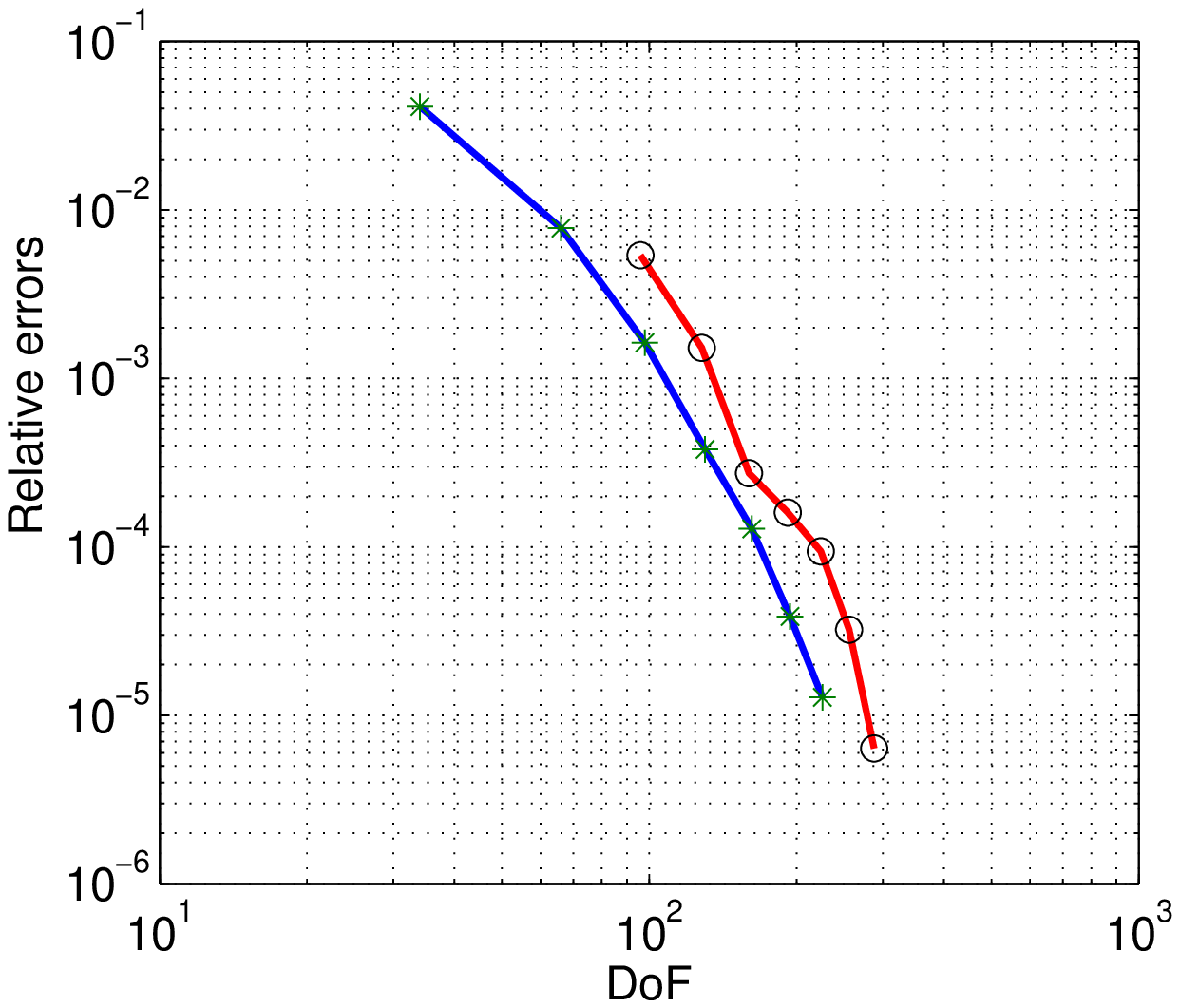}}
\caption{Example 2: relative errors of Table \ref{ex2_tab3}  vs
DoF with $h=1/32$ ('$*$' B-splines, '$\circ$' Lagrangian basis). }
\label{errorivsdof_IIraff}
\end{figure}

Then, in Tables \ref{ex2_tab2} and \ref{ex2_tab3}, results
obtained considering two successive halving of the mesh size $h$
are reported. The corresponding error behaviors are shown in
Figures~\ref{errorivsdof_Iraff} and \ref{errorivsdof_IIraff},
respectively:
the rate of convergence of both methods are almost equal, although, using IGA-SGBEM we can achieve the same error with fewer DoFs w.r.t. C-SGBEM. This is due to the fact that in the so-called \textit{$p$-version} of the
Galerkin BEM (\cite{suri, postell}), where the accuracy is reached
fixing the mesh  and elevating the  degree of the piecewise
polynomial basis, the error mainly depends on this degree
(order).

\begin{figure}
\centerline{\includegraphics[width=2.7in]{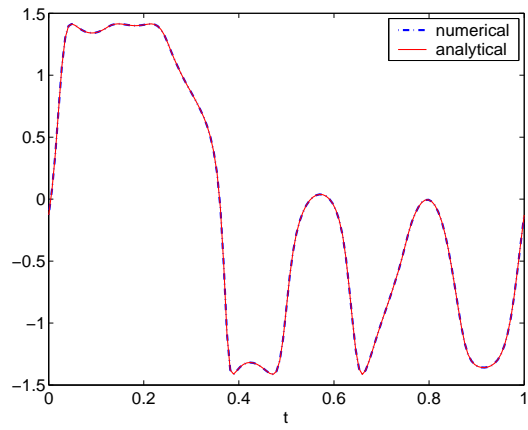}}
\caption{Example 2: approximate solution obtained with the
B-spline basis of degree $9$ and $h=1/16$, together with the
analytical solution.} \label{sol_free_form}
\end{figure}

Further, note that in these simulations the conditioning of B-splines systems are worse
than the corresponding Lagrangian ones. Anyway, the remarkable shape reproduction capability of
our scheme is underlined in Figure \ref{sol_free_form} which shows
the approximate solution obtained with $h=1/16$ and the B-spline
basis of degree $9$ together with the analytical solution.

Finally, in order to put in evidence possible benefits of our
approach, we fix now the degree of the piecewise polynomial spaces
equal to 3 and compare the results obtained using $C^2$ B-spline
basis in IGA-SGBEM and $C^0$ Lagrangian basis in C-SGBEM. The
comparison is first done in Table \ref{ex2_tab4}, with respect
to $h$, where results obtained using $C^1$ B-spline basis in
IGA-SGBEM are also given. We note that the errors with the
IGA-SGBEM approach are slightly worse, but the degrees of freedom
are remarkably lower.

\begin{figure}[bht]
\centerline{\includegraphics[width=2.8in]{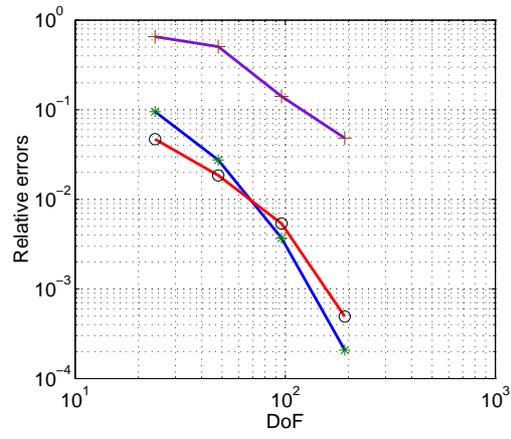}}
\caption{Example 2: relative errors of Table
\ref{ex2_tab5} ('$*$' IGA-SGBEM, '$\circ$' C-SGBEM, '$+$' S-SGBEM).} \label{errorivsdof_parita_dof}
\end{figure}

The benefits of our approach can be better
observed looking at Table \ref{ex2_tab5}, where indeed the
comparison is done with respect to DoF. Such a comparison can be
achieved by selecting, for a given value of DoF, a suitable mesh
size $h$. To complete this benchmark, on the right of Table \ref{ex2_tab5}, for the same
DoF, relative errors obtained using $L^2$ cubic Lagrangian basis
on piecewise linear approximation $\tilde{\Gamma}_h$ of the
boundary $\Gamma$ are reported.
All the errors of these last Table are plotted in Figure
\ref{errorivsdof_parita_dof} with respect to DoF. This graphic
reveals the inferiority of the standard SGBEM approach
with respect to IGA and curvilinear SGBEMs, due to the introduced
approximation of the boundary.

\medskip

\noindent {\bf Example 3.}

\begin{figure}
\centerline{\includegraphics[width=3.in]{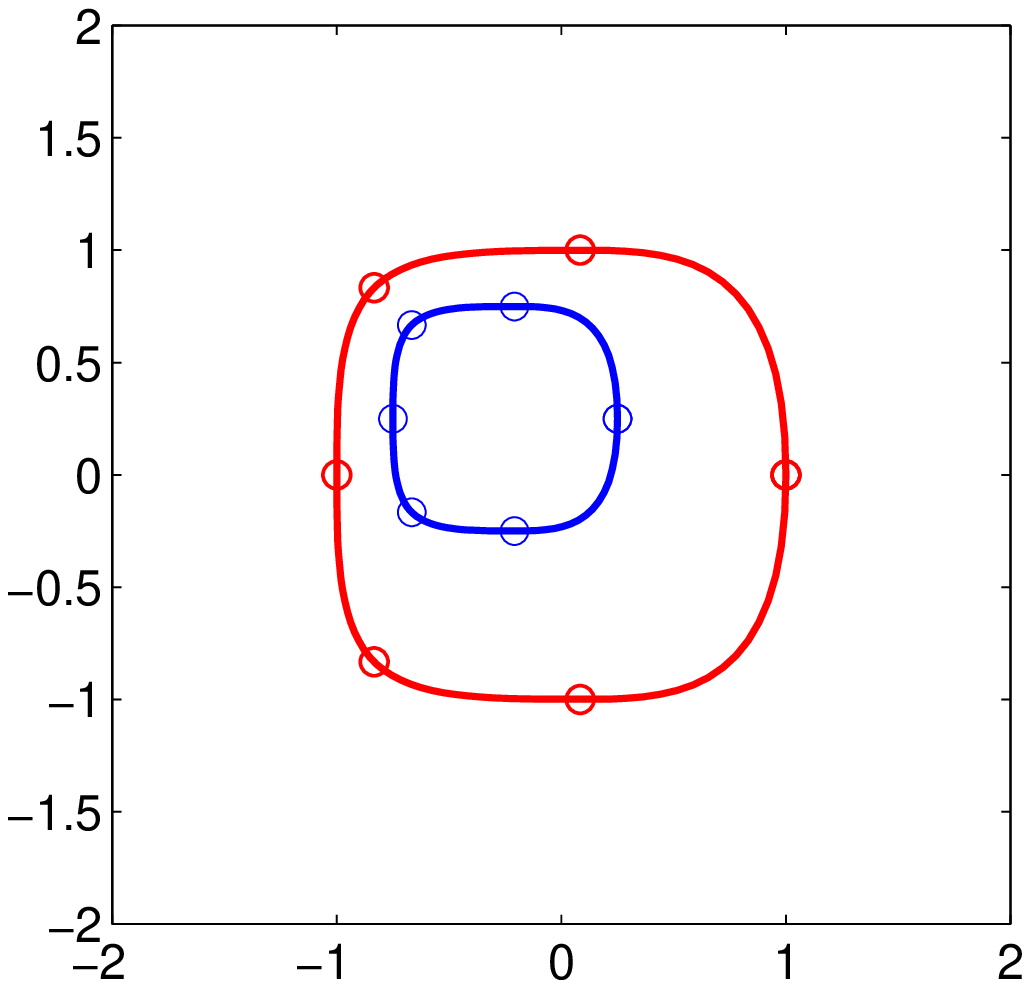}
\includegraphics[width=3in]{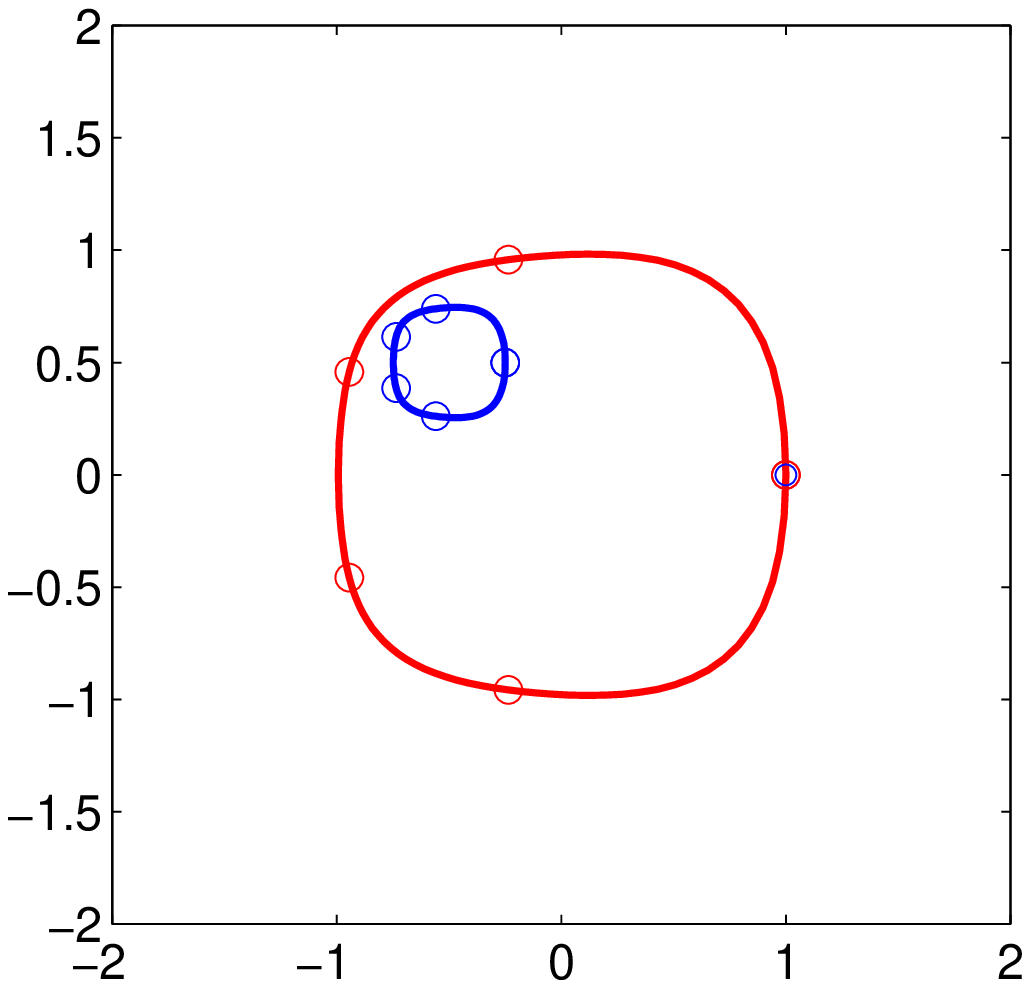}}
\caption{Example 3: Domains with holes. The two closed boundary
curves on the left ($A$) are cubics and those on the right ($B$)
are quartics. The nodal (mesh) points of the spline curves are
marked with the symbol '$\circ$'.} \label{D_domain_holes}
\end{figure}

Since one of the major strengths of BEM  approach (with respect to
FEM) is its ability of easily treating domains with holes, let us
now consider the two trimmed domains depicted in Figure
\ref{D_domain_holes}. For the domain on the left ($A$), the two
boundary curves are represented by cubic B-splines with extended
knot vector
$$ T_4=[0,\, 0,\, 0,\, 0,\, 1/6,\,
2/6,\, 3/6,\, 4/6,\, 5/6,\, 1,\, 1,\, 1,\, 1 ],
$$
while the curves defining the domain on the right ($B$), are
quartic B-splines with extended
knot vector
$$T_5=[0,\, 0,\, 0,\, 0,\, 0,\, 1/5,\, 2/5,\, 3/5,\,
4/5,\, 1,\, 1,\, 1,\, 1,\, 1].$$
The coordinates of the control points associated to the external
and internal boundary curves are collected in the following
matrices:
$$
Q_{\rm est_A}=\Big[\begin{array}{c c c c c c c c c }
 1 &    1   &  0   & -1   & -1  &  -1  &   0   &  1   &  1\\
     0   &  1    & 1  &   1  &   0  &  -1  &  -1  &  -1  &   0
\end{array}\Big]\,,
$$
$$
Q_{\rm int_A}=\Big[\begin{array}{c c c c c c c c c }
 0.25 &    0.25   &  -0.25   & -0.75   & -0.75  &  -0.75  &   -0.25   &  0.25   &  0.25\\
     0.25   &  -0.25    & -0.25  &   -0.25  &   0.25  & 0.75  &  0.75  &  0.75  &   0.25
\end{array}\Big]\,.
$$
$$
Q_{\rm est_B}=\Big[\begin{array}{c c c c c c c c c }
 1 &    1   &  0   & -1   & -1  &  -1  &   0   &  1   &  1\\
     0   &  1    & 1  &   1  &   0  &  -1  &  -1  &  -1  &   0
\end{array}\Big]\,,
$$
$$
Q_{\rm int_B}=\Big[\begin{array}{c c c c c c c c c }
 -0.25 &    -0.25   &  -0.5   & -0.75   & -0.75  &  -0.75  &   -0.5   &  -0.25   &  -0.25\\
     0.5   &  0.25    & 0.25  &   0.25  &   0.5  & 0.75  &  0.75  &  0.75  &   0.5
\end{array}\Big]\,.
$$

For both domains a mixed BVP is considered, where a Dirichlet
condition $u^*=1$ is assigned on the interior boundary, while a
Neumann condition $q^*=0$ is prescribed on the exterior boundary.
This configuration can model a stationary heat conduction problem,
where a constant temperature on the inner wall and a zero heat
flux on the outer wall are given.

Focusing on domain $A$, we have tested our IGA-SGBEM approach just
using the cubic spline space used to define the boundary curves
($h=1/6$). The resulting linear system is of order $16$ and the
approach produces an approximate solution with an absolute error $E_M$ in maximum
norm equal to $1.9451\, 10^{-5}$ for what concerns the recovered
flux $q$ and equal to $2.1153\, 10^{-5}$ for what concerns the
recovered potential $u$. If we use cubic C-SGBEM, instead, we have
to solve a liner system of order $36$ to reproduce the same error
order.

Regarding domain $B$, we have again tested our scheme just
considering the quartic spline space used to describe the
boundaries ($h=1/5$). Again we end up with a linear system of
order $16$, which produces an approximate solution with an
absolute error $E_M$ in maximum norm equal to $3.7748\, 10^{-5}$ for
what concerns the recovered flux $q$ and equal to $3.7301\,
10^{-5}$ for what concerns the recovered potential $u$. If we use
quartic C-SGBEM, instead, we have to solve a liner system of order
$40$ to reproduce the same error order.

\vspace{0.1in} \noindent \textit{Remark}. Here, we have chosen to evaluate the absolute error $E_M$ in maximum norm
instead of \eqref{E} since the considered mixed
boundary conditions allow the BVP to have the constant solution
$u=1$. The
obtained errors are due to the approximation of weakly
singular, singular and hypersingular double integrals by
means of the already mentioned quadrature formulas (\cite{Aimi1}).

\medskip


\noindent {\bf Example 4.}

Let us conclude this Section, considering a Dirichlet BVP for the
Laplace equation exterior to the arc of parabola $\Gamma=\{{\bf
x}=(x_1,x_2)|\,  x_1=t,\,x_2=1-t^2, \,t\in [a,b]=[-1,1]$\},
representable by means of quadratic B-splines related to the
extended knot vector
 $$T_6=\begin{array}{c c c c c c}
[-1 &-1& -1& 1& 1& 1]
\end{array}$$
and to the control points $\Q_i,\, i=0,\cdots,2$, whose
coordinates are collected in the following matrix:
 $$ Q=\Big[\begin{array}{c c c}
-1 &0& 1\\0& 2& 0
\end{array}\Big].$$

The considered differential problem can model the electrostatic
problem of finding the electric potential around a condenser,
whose two faces are so near one another to be considered as
overlapped, knowing the electric potential only on the condenser.
Here the Dirichlet datum is given in such a way that the solution
of the related boundary integral equation is explicitly known and
reads $[q(\x)]_\Gamma=\sqrt{1+4 x_1^2}$.

\begin{figure}
\centerline{\includegraphics[width=2.8in]{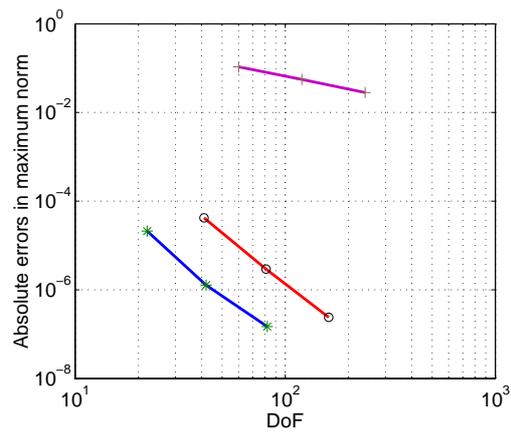}}
\caption{Example 4: Absolute errors, related to quadratic basis
functions, for different value of $h$ ('$*$' $C^1$ B-splines,
'$\circ$' $C^0$ Lagrangian basis, '$+$' $L^2$ Lagrangian basis on
$\tilde{\Gamma}_h$).} \label{errorevsdof_IIesempio_II}
\end{figure}

The comparison reported in Table \ref{tabella4bis}, for different values of the parameter $h$, which uniformly decomposes the parameter interval
$[-1,1]$, involves $C^1$ quadratic B-spline basis for IGA-SGBEM,
$C^0$ quadratic Lagrangian basis for C-SGBEM and $L^2$  quadratic
Lagrangian basis on piecewise linear approximation
$\tilde{\Gamma}_h$ of $\Gamma$ for S-SGBEM. Together with DoF and
spectral condition numbers of the associated matrices, we show the
absolute errors $E_M$ in maximum norm. These errors are visualized
in Figure \ref{errorevsdof_IIesempio_II} with respect to DoF.

At last, for this example, in Table \ref{ultima} we show a comparison between
Galerkin IGA-BEM described in this paper and collocation IGA-BEM,
where collocation is done at the
Greville abscissae as in \cite{Simpson}. It turns out that, for a fixed
discretization parameter $h$, the Galerkin technique is more
accurate than the collocation one, while the matrix condition
number of the latter is better, even if the symmetry property
useful in the coupling with FEM (\cite{Wendland2, Costabel0}) is
lost. Both techniques, as shown in this Table, satisfy the
estimates given in \cite{Rannacher} for what concerns the decay of
Galerkin error $E_M$, which, for smooth boundaries and
sufficiently regular data, behaves as $O(h^k)$, being $k$ the order of the fixed B-spline basis.

\section{Conclusions}
In this work we studied from a numerical point of view an
Isogeometric Symmetric Galerkin Boundary Element Method, which we
called IGA-SGBEM, dealing with the reference 2D Laplace problem,
on domains having different shapes. In particular our aim was to
compare the performances of such approach not only with those of
standard SGBEM (where the boundary of the domain is approximated
by polygonal lines), but also with those of a more advanced SGBEM,
namely curvilinear SGBEM, which is capable of retaining the exact
boundary.

The potential strength and superiority of the presented approach
has been confirmed by all the numerical tests, where smooth and
non smooth interior domains as well as domains with holes or
unbounded domains exterior to an open arc have been considered.
The only drawback of IGA-SGBEM is, in few cases involving very
long boundary elements, the worse conditioning of the
discretization linear system matrix, which is probably due to the
larger support of B-splines with respect to that one of the
Lagrangian basis functions.

In order to better exploit the potentiality of the isogeometric
approach combined with SGBEM we plan, as a future work, to
extend the analysis to non polynomial spline spaces, able to
represent exactly complex shapes. This can be achieved by
considering generalized B-splines (\cite{MPS'11}). The
extension of our approach to 3D problems would constitute a
further challenge, where its major appeal with respect to
classical IGA-FEM could be more evident.

\section*{Acknowledgements}
This work has been partially supported by INdAM, through GNCS
research projects.

\newpage
\begin{table} [htp]
\begin{center}
\begin{tabular} {|c|c c c|}
\hline
 $h$ & DoF & ${\rm cond.}$ & $E$ \\
\hline
 $1$ & $13$ & $2.53\,10^{2}$& $3.32\,10^{-1}$  \\
 $1/2$ & $22$ &$3.05\,10^{2}$ & $1.61\,10^{-1}$ \\
 $1/4$ & $40$ & $5.82\,10^{2}$& $1.08\,10^{-1}$ \\
 $1/8$ & $76$ & $1.23\,10^{3}$& $7.62\,10^{-2}$ \\
\hline
\end{tabular}
\end{center}
\caption{Example 1: results obtained by quadratic B-splines
starting from extended knot vector $T_1$.} \label{ex1_tab1}
\end{table}
\begin{table} [htp]
\begin{center}
\begin{tabular} {|c|c c c| c c c |}
\hline
 & &  {\rm IGA-SGBEM} & & &  {\rm C-SGBEM}  & \\
\hline
$h$ & DoF & ${\rm cond.}$ & $E$ & DoF & ${\rm cond.}$ & $E$\\
\hline
$1$ &  $15$ & $4.23\,10^{2}$ & $2.06\,10^{-1}$ & $21$ & $1.68\,10^{2}$ &$4.56\,10^{-2}$ \\
 $1/2$ &  $24$ & $3.55\,10^{2}$ & $5.20\,10^{-2}$ & $39$ & $2.68\,10^{2}$ & $1.84\,10^{-2}$ \\
 $1/4$ &  $42$ & $5.68\,10^{2}$ & $1.83\,10^{-2}$  & $75$ &$5.38\,10^{2}$ & $2.35\,10^{-2}$ \\
 $1/8$ &  $78$ & $1.21\,10^{3}$ & $1.84\,10^{-2}$  & $147$ &$1.09\,10^{3}$ & $3.34\,10^{-2}$ \\
\hline
\end{tabular}
\caption{Example 1: comparison between results obtained with
quadratic $L^2(0,9)\cap C^1[1,8]$ B-splines/ $L^2(0,9)\cap
C^0[1,8]$ Lagrangian basis, varying $h$.} \label{ex1_tab2}
\end{center}
\end{table}
\begin{table} [htp]
\begin{center}
\begin{tabular} {|c|c c c|c c c|}
\hline
$h=1/8$& & IGA-SGBEM & & & C-SGBEM  & \\
\hline
 ${\rm degree}$ & DoF& ${\rm cond.}$ & $E$ & DoF&
 ${\rm cond.}$ & $E$\\
\hline
 $3$ & $10$ & $6.61\,10^2$& $3.37\,10^{-1}$ & $24$ & $5.47\,10^2$& $4.69\,10^{-2}$ \\
 $4$ & $18$ &$4.07\,10^3 $ & $1.26\,10^{-1}$ & $32$ &$1.21\,10^3$ & $3.43\,10^{-2}$\\
 $5$ & $26$ & $1.89\,10^4$& $6.55\,10^{-2}$ & $40$ & $2.02\,10^3$& $2.34\,10^{-2}$\\
 $6$ & $34$ & $9.19\,10^4$& $3.23\,10^{-2}$ & $48$ & $4.42\,10^3$& $1.56\,10^{-2}$ \\
 $7$ & $42$ &$4.40\,10^5 $ & $1.71\,10^{-2}$ & $56$ &$5.50\,10^3$ & $1.17\,10^{-2}$\\
 $8$ & $50$ & $2.08\,10^6$& $1.03\,10^{-2}$ & $64$ & $1.88\,10^4$& $8.51\,10^{-3}$\\
 $9$ & $58$ & $9.84\,10^6$& $3.57\,10^{-3}$ & $72$ & $1.45\,10^4$& $3.26\,10^{-3}$\\
\hline
\end{tabular}
\end{center}
\caption{Example 2: comparison between results obtained with
IGA-SGBEM based on $C^2$ B-splines and C-SGBEM based on $C^0$
Lagrangian basis, for different degrees of the piecewise
polynomial basis and $h=1/8$.} \label{ex2_tab1}
\end{table}
\begin{table} [htp]
\begin{center}
\begin{tabular} {|c|c c c|c c c|}
\hline
$h=1/1$6& & IGA-SGBEM& & & C-SGBEM & \\
\hline
 ${\rm degree}$ & DoF& ${\rm cond.}$ & $E$
 & DoF& ${\rm cond.}$ & $E$\\
\hline
 $3$ & $18$ & $1.18\,10^3$& $1.28\,10^{-1}$ & $48$ & $1.49\,10^3$& $1.85\,10^{-2}$ \\
 $4$ & $34$ &$7.78\,10^3 $ & $4.06\,10^{-2}$ & $64$ &$3.36\,10^3$ & $1.13\,10^{-2}$\\
 $5$ & $50$ & $3.72\,10^4$& $1.44\,10^{-2}$ & $80$ & $5.30\,10^3$& $6.32\,10^{-3}$\\
 $6$ & $66$ & $1.90\,10^5$& $4.92\,10^{-3}$ & $96$ & $1.16\,10^4$& $2.52\,10^{-3}$ \\
 $7$ & $82$ &$9.69\,10^5 $ & $1.97\,10^{-3}$ & $112$ &$1.35\,10^4$ & $1.10\,10^{-3}$\\
 $8$ & $98$ & $4.98\,10^6$& $7.31\,10^{-4}$ & $128$ & $4.93\,10^4$& $4.31\,10^{-4}$\\
 $9$ & $114$ & $2.51\,10^7$& $4.08\,10^{-4}$ & $144$ & $3.49\,10^4$& $3.74\,10^{-4}$\\
\hline
\end{tabular}
\end{center} \caption{Example 2: comparison between results obtained with
IGA-SGBEM based on $C^2$ B-splines and C-SGBEM based on $C^0$
Lagrangian basis, for different degrees of the piecewise
polynomial basis and $h=1/16$.} \label{ex2_tab2}
\end{table}
\begin{table} [htp]
\begin{center}
\begin{tabular} {|c|c c c|c c c|}
\hline
$h=1/32$& & IGA-SGBEM& & &  C-SGBEM& \\
\hline
 ${\rm degree}$ & DoF& ${\rm cond.}$ & $E$ & DoF&
 ${\rm cond.}$ & $E$\\
\hline
 $3$ &$34$&$2.44\,10^3$&$4.11\,10^{-2}$&$96$ &$4.21\,10^3$&$5.38\,10^{-3}$ \\
 $4$ &$66$&$2.21\,10^4$ & $7.81\,10^{-3}$ & $128$ &$8.84\,10^3$ & $1.52\,10^{-3}$\\
 $5$ &$98$&$1.13\,10^5$ & $1.63\,10^{-3}$ & $160$ & $1.32\,10^4$& $2.74\,10^{-4}$\\
 $6$ &$130$&$5.90\,10^5$& $3.80\,10^{-4}$ & $192$ & $2.77\,10^4$& $1.60\,10^{-4}$ \\
 $7$ &$162$&$3.08\,10^6$ & $1.28\,10^{-4}$ & $224$ &$3.09\,10^4$ & $9.40\,10^{-5}$\\
 $8$ &$194$&$1.59\,10^7$& $3.86\,10^{-5}$ & $256$ & $1.12\,10^5$& $3.23\,10^{-5}$\\
 $9$ &$226$&$8.08\,10^7$& $1.28\,10^{-5}$ & $288$ & $7.53\,10^4$& $6.36\,10^{-6}$\\
\hline
\end{tabular}
\end{center} \caption{Example 2: comparison between results obtained with
IGA-SGBEM based on $C^2$ B-splines and C-SGBEM based on $C^0$
Lagrangian basis, for different degrees of the piecewise
polynomial basis and $h=1/32$.} \label{ex2_tab3}
\end{table}
\begin{table} [htp]
\hspace{-0.25in}
\begin{tabular} {|c|c c c|c c c|c c c|}
\hline
& & $C^2$ IGA-SGBEM& & & $C^1$ IGA-SGBEM & & & $C^0$ C-SGBEM  & \\
\hline
 $h$ & DoF& ${\rm cond.}$ & $E$ & DoF& ${\rm cond.}$
 & $E$ &DoF& ${\rm cond.}$ & $E$\\
\hline
 $1/8$ & $10$ & $6.61\,10^2$& $3.37\,10^{-1}$ & $17$ & $1.29\,10^3$&
$1.30\,10^{-1}$& $24$ & $5.47\,10^2$& $4.69\,10^{-2}$ \\
 $1/16$ & $18$ &$1.18\,10^3 $ & $1.28\,10^{-1}$ & $25$ &$2.17\,10^3 $ &
$6.23\,10^{-2}$&$48$ &$1.49\,10^3$ & $1.85\,10^{-2}$\\
 $1/32$ & $34$ & $2.44\,10^3$& $4.11\,10^{-2}$ & $41$ & $3.99\,10^3$&
$1.95\,10^{-2}$&  $96$ & $4.21\,10^3$& $5.38\,10^{-3}$\\
 $1/64$ & $66$ & $7.48\,10^3$& $8.29\,10^{-3}$ & $73$ & $1.19\,10^4$&
$2.63\,10^{-3}$& $192$ & $1.10\,10^4$& $4.92\,10^{-4}$ \\
\hline
\end{tabular}
\caption{Example 2: comparison among results obtained with cubic
$C^2$ B-splines (left), cubic $C^1$ B-splines (middle) and cubic
$C^0$ Lagrangian basis (right).} \label{ex2_tab4}
\end{table}
\begin{table} [htp]
\begin{center}
\begin{tabular} {|c|c c c|c c c|c c c|}
\hline
& & IGA-SGBEM & & & C-SGBEM & & & S-SGBEM & \\
\hline
 DoF & $h$& ${\rm cond.}$ & $E$ & $h$& ${\rm cond.}$
 & $E$ & $h$& ${\rm cond.}$ & $E$\\
\hline
 $24$ & $1/22$ & $1.63\,10^3$& $9.53\,10^{-2}$ & $1/8$ & $5.47\,10^2$& $4.69\,10^{-2}$ & $1/6$ & $3.01\,10^2$& $6.54\,10^{-1}$\\
 $48$ & $1/46$ &$4.25\,10^3 $ & $2.75\,10^{-2}$ & $1/16$ &$1.49\,10^3$ & $1.85\,10^{-2}$ & $1/12$ &$1.26\,10^3 $ & $5.06\,10^{-1}$\\
 $96$ & $1/94$ & $1.46\,10^4$& $3.70\,10^{-3}$ & $1/32$ & $4.21\,10^3$& $5.38\,10^{-3}$ & $1/24$ & $1.59\,10^3$& $1.41\,10^{-1}$ \\
 $192$ & $1/190$ & $4.57\,10^4$& $2.08\,10^{-4}$ & $1/64$ & $1.10\,10^4$& $4.92\,10^{-4}$ & $1/48$ & $5.29\,10^3$& $4.82\,10^{-2}$\\
\hline
\end{tabular}
\end{center} \caption{Example 2: comparison between results obtained with cubic $C^2$ B-splines (left), $C^0$
Lagrangian basis (middle), $L^2$
Lagrangian basis on ${\tilde \Gamma}_h$ (right), for different values of DoF.} \label{ex2_tab5}
\end{table}
\begin{table} [htp]
\hspace{0.1in}
\begin{tabular} {|c|c c c| c c c |c c c|}
\hline
 &&IGA-SGBEM&&&C-SGBEM&&&S-SGBEM&\\
\hline
$h$ & DoF&  cond. & $E_M$ & DoF& cond. & $E_M$& DoF& cond & $E_M$ \\
\hline
 $1/10$ &  $22$ & $1.87\,10^{2}$ & $2.11\,10^{-5}$ & $41$ & $2.33\,10^{2}$ &$4.20\,10^{-5}$ & $60$ &$1.01\,10^{3}$&$1.07\,10^{-1}$\\
 $1/20$ &  $42$ & $4.57\,10^{2}$ & $1.27\,10^{-6}$ & $81$ & $5.00\,10^{2}$ & $2.92\,10^{-6}$ & $120$& $2.09\,10^{3}$&$5.55\,10^{-2}$\\
 $1/40$ &  $82$ & $1.01\,10^{3}$ & $1.48\,10^{-7}$  & $161$ &$1.04\,10^{3}$ & $2.40\,10^{-7}$ & $240$  & $4.27\,10^{3}$ &$2.82\,10^{-2}$\\
\hline
\end{tabular}
\caption{Example 4: comparison between IGA-SGBEM, C-SGBEM, S-SGBEM, based on quadratic piecewise polynomial basis functions, for different values of $h$.}
\label{tabella4bis}
\end{table}
\begin{table} [htp]
\begin{center}
\begin{tabular} {|c c| c c c| c c c|}
\hline
 & & & Galerkin IGA-BEM  & & & collocation IGA-BEM  & \\
\hline
$h$ & DoF & ${\rm cond.}$ & $E_M(h)$ & $\log_2\Big (\frac{E_M(2h)}{E_M(h)}\Big )$ & ${\rm cond.}$ & $E_M(h)$ &$\log_2\Big (\frac{E_M(2h)}{E_M(h)}\Big )$\\
\hline
 $1/5$ &  $12$ & $7.19\,10^{1}$ & $4.03\,10^{-4}$ & $-$ & $1.73\,10^{1}$ & $4.88\,10^{-4}$ & $-$\\
 $1/10$ &  $22$ & $1.87\,10^{2}$ & $2.11\,10^{-5}$ & $4.26$ & $3.63\,10^{1}$ & $5.74\,10^{-5}$ & $3.09$\\
 $1/20$ &  $42$ & $4.57\,10^{2}$ & $1.27\,10^{-6}$  & $4.05$ &$7.70\,10^{1}$ & $6.94\,10^{-6}$ & $3.05$\\
 $1/40$ &  $82$ & $1.01\,10^{3}$ & $1.48\,10^{-7}$  & $3.10$ &$1.58\,10^{2}$ & $8.58\,10^{-7}$ & $3.01$\\
 $1/80$ &  $162$ & $2.12\,10^{3}$ & $1.81\,10^{-8}$  & $3.03$ &$3.20\,10^{2}$ & $1.07\,10^{-7}$ & $3.00$\\
\hline
\end{tabular}
\caption{Example 4: comparison between Galerkin and collocation
IGA-BEM results, varying $h$.} \label{ultima}
\end{center}
\end{table}
\clearpage

\end{document}